\documentclass[12pt]{article}
\usepackage{makeidx}
\usepackage{amssymb, amsmath, amsfonts, enumerate}
\usepackage{setspace, graphicx, caption, subcaption, float, relsize, rotating, color, verbatim, multirow}
\usepackage[round, comma, longnamesfirst]{natbib}
\usepackage[compact]{titlesec}

\makeatletter
\def\maketag@@@#1{\hbox{\m@th\normalfont\normalsize#1}}
\makeatother
\def\maketag@@@#1{\hbox{\m@th\normalfont\normalsize#1}}
\makeatother

\allowdisplaybreaks
\newtheorem {theorem}{Theorem}[section]
\newtheorem {assumption}{Assumption}

\newtheorem{lemma}[theorem]{Lemma}

\newtheorem{remark}{Remark}

\input{tcilatex}

\begin{document}

\title{Asymptotic Theory for the Maximum of an Increasing Sequence of
Parametric Functions}
\author{Jonathan B. Hill \thanks{%
Department of Economics, University of North Carolina at Chapel Hill;
jbhill@email.unc.edu.} \\
%EndAName
University of North Carolina}
\date{\today \\
}
\maketitle

\begin{abstract}
\cite{HillMotegi2017} present a new general asymptotic theory for the
maximum of a random array $\{\mathcal{X}_{n}(i)$ $:$ $1$ $\leq $ $i$ $\leq $
$\mathcal{L}\}_{n\geq 1}$, where each $\mathcal{X}_{n}(i)$ is assumed to
converge in probability as $n$ $\rightarrow $ $\infty $. The array dimension
$\mathcal{L}$ is allowed to increase with the sample size $n$. Existing
extreme value theory arguments focus on observed data $\mathcal{X}_{n}(i)$,
and require a well defined limit law for $\max_{1\leq i\leq \mathcal{L}}|\mathcal{X}_{n}(i)|$ by restricting dependence across $i$. The high
dimensional central limit theory literature presumes approximability by a
Gaussian law, and also restricts attention to observed data. \cite{HillMotegi2017} do not require $\max_{1\leq i\leq \mathcal{L}_{n}}|\mathcal{X}_{n}(i)|$ to have a well defined limit nor be approximable by a Gaussian
random variable, and we do not make any assumptions about dependence across $i$. We apply the theory to filtered data when the variable of interest $\mathcal{X}_{n}(i,\theta _{0})$ is not observed, but its sample counterpart $\mathcal{X}_{n}(i,\hat{\theta}_{n})$ is observed where $\hat{\theta}_{n}$
estimates $\theta _{0}$. The main results are illustrated by looking at unit
root tests for a high dimensional random variable, and a residuals white
noise test.\medskip \newline
\textbf{Keywords} : maximum of multivariate mean, high dimensionality,
non-Gaussian approximation.\medskip \newline
\textbf{MSC2010 subject classifications} : 62E17, 62F40, 62M10, 05D10.
\end{abstract}

\setstretch{1.6}

\section{Introduction\label{sec:intro}}

Consider an array of random variables on a probability measure space $%
(\Omega ,\mathcal{F},\mathcal{P})$:
\begin{equation*}
\left\{ \mathcal{X}_{n}(i),\mathcal{Y}_{n}(i):1\leq i\leq \mathcal{I}%
_{n}\right\} _{n\geq 1},
\end{equation*}%
where $\{\mathcal{I}_{n}\}_{n\geq 1}$ is a sequence of positive integers, $%
\mathcal{I}_{n}$ $\rightarrow $ $\infty $ as $n$ $\rightarrow $ $\infty $.
Under the assumption that $\mathcal{X}_{n}(i)$ $\overset{p}{\rightarrow }$ $0
$ as $n$ $\rightarrow $ $\infty $\ for each $i$, \cite{HillMotegi2017} prove
for some sequence of positive integers $\{\mathcal{L}_{n}\}$ that satisfies $%
\mathcal{L}_{n}$ $\rightarrow $ $\infty $:%
\begin{equation}
\max_{1\leq i\leq \mathcal{L}_{n}}\left\vert \mathcal{X}_{n}(i)\right\vert
\overset{p}{\rightarrow }0.  \label{max|X|}
\end{equation}%
By operating on (\ref{max|X|}), for any two arrays of $\mathcal{F}$%
-measurable random variables $\{\mathcal{X}_{n}(i),\mathcal{Y}_{n}(i)$ $:$ $1
$ $\leq $ $i$ $\leq $ $\mathcal{I}_{n}\}_{n\geq 1}$, if $\mathcal{X}_{n}(i)$
$-$ $\mathcal{Y}_{n}(i)$ $\overset{p}{\rightarrow }$ $0$ for each $i$ \cite%
{HillMotegi2017} then show that the difference in maxima satisfy:%
\begin{equation}
\left\vert \max_{1\leq i\leq \mathcal{L}_{n}}\left\vert \mathcal{X}%
_{n}(i)\right\vert -\max_{1\leq i\leq \mathcal{L}_{n}}\left\vert \mathcal{Y}%
_{n}(i)\right\vert \right\vert \overset{p}{\rightarrow }0.  \label{max|X-Y|}
\end{equation}%
Examples are a sample mean $\mathcal{X}_{n}(i)$ $=$ $1/n%
\sum_{t=1}^{n}x_{t}(i)$ or standardized means $\mathcal{X}_{n}(i)$ $=$ $1/%
\sqrt{n}\sum_{t=1}^{n}x_{t}(i)$ and $\mathcal{Y}_{n}(i)$ $=$ $1/\sqrt{n}%
\sum_{t=1}^{n}y_{t}(i)$, where $\{x_{t}(i),y_{t}(i)\}_{t=1}^{n}$ are the
sample paths of some processes on $(\Omega ,\mathcal{F},\mathcal{P})$. This
has been studied extensively in the Gaussian approximation and high
dimensional Gaussian central limit theory literatures, discussed below. At
the highest level of generality we work with maxima of positive values in
order to exploit convenient inequalities.

The focus of the present paper is to extend the key ideas of \cite%
{HillMotegi2017} to filtered residuals, and to apply the result to a unit
root test and a white noise test. The authors do not impose any restrictions
on dependence in $(\mathcal{X}_{n}(i),\mathcal{Y}_{n}(i))$ across
coordinates $i$, nor do we require $(\mathcal{X}_{n}(i),\mathcal{Y}_{n}(i))$
to belong to a specific domain of attraction. For example, in the normalized
mean case we neither require either $(\mathcal{X}_{n}(i),\mathcal{Y}_{n}(i))$
to be Gaussian nor converge in law to Gaussian random variables. The
generality arises from a new result for convergence of high dimensional
arrays which does not rely on probabilistic properties, although naturally
lends itself to probabilistic applications.

The approach differs from standard weak convergence methods when applied to
an array of standardized means $\{\mathcal{X}_{n}(i)$ $:$ $1$ $\leq $ $i$ $%
\leq $ $\mathcal{I}_{n}\}_{n\geq 1}$, e.g. $\mathcal{X}_{n}(i)$ $=$ $1/\sqrt{%
n}\sum_{t=1}^{n}x_{t}(i)$. Weak convergence for $\{\mathcal{X}_{n}(i)$ $:$ $i
$ $\in $ $N\}$ in the broad sense of \citet{HoffJorg1984,HoffJorg1991} is a
potential option. It is known that such weak convergence to a Gaussian
limit, with a versin that has uniformly bounded and uniformly continuous
sample paths, is equivelant to pointwise convergence and the existence of a
pseudo metric $d$ on $N$\ such that $(N,d)$ is a totally bounded pseudo
metric space and a stochastic equicontinuity property based on $d$ holds.
See \citet{Dudley1978,Dudley1984} and \citet[Chapters 9-10]{Pollard1990}. If
$d$ is the Euclidean distance, for example, then $(N,d)$ is not totally
bounded. \cite{HillMotegi2017} take a different approach that completely
sidesteps the approach of \citet{HoffJorg1984,HoffJorg1991} by first
demonstrating pointwise convergence $\mathcal{X}_{n}(i)$ $\overset{d}{%
\rightarrow }$ $\mathcal{X}(i)$. They then work with a probability
construction for the array $\{\mathcal{X}_{n}(i)$ $-$ $\mathcal{X}(i)$ $:$ $1
$ $\leq $ $i$ $\leq $ $\mathcal{I}_{n}\}_{n\geq 1}$ and apply our general
result for convergence of arrays to be able to show $|\max_{1\leq i\leq
\mathcal{L}_{n}}|\mathcal{X}_{n}(i)|$ $-$ $\max_{1\leq i\leq \mathcal{L}%
_{n}}|\mathcal{X}(i)||$ $\overset{p}{\rightarrow }$ $0$.

The primary tools used here date in some form to seminal theory developed by
\cite{Ramsey1930} and its implications for monotone subsequences and
convergence, cf. the Erd\H{o}s-Szekeres theorem \citep{ErdosSzekeres1935}.
\cite{HillMotegi2017} significantly augment a convergence result for
non-stochastic arrays presented in \citet[Lemma
1]{BoehmeRosenfeld1974} in key ways. The latter claim that if an array $\{%
\mathcal{A}_{k,n}$ $:$ $1$ $\leq $ $k$ $\leq $ $\mathcal{I}_{n}\}_{n\geq 1}$%
, where $\mathcal{I}_{n}$ $\rightarrow $ $\infty $ as $n$ $\rightarrow $ $%
\infty $, lies in a first countable topological space, and $%
\lim_{k\rightarrow \infty }\lim_{n\rightarrow \infty }\mathcal{A}_{k,n}$ $=$
$0$, then $\lim_{l\rightarrow \infty }\mathcal{A}_{\mathcal{L}(n_{l}),n_{l}}$
$=$ $0$ for some infinite subsequence $\{n_{l}\}_{l=1}^{\infty }$ of
positive integers, and some mapping $\mathcal{L}(n_{l})$ $\rightarrow $ $%
\infty $. (Recall that any metric space is a first countable topological
space.) In order to use the result for maxima, we extend $%
\{n_{l}\}_{l=1}^{\infty }$ to $\mathbb{N}$ and therefore achieve $%
\lim_{n\rightarrow \infty }\mathcal{A}_{\mathcal{L}(n),n}$ $=$ $0$ in Lemma %
\ref{lm:array_conv}. We require additional assumptions that lends itself to
deriving (\ref{max|X|}) and (\ref{max|X-Y|}). A practical application
includes when $\mathcal{A}_{k,n}$ $=$ $\int_{0}^{1}P(1$ $-$ $e^{-\max_{1\leq
i\leq k}|\mathcal{X}_{n}(i)|}$ $>$ $\epsilon )d\epsilon $, the foundation
for showing $\max_{1\leq i\leq \mathcal{L}_{n}}|\mathcal{X}_{n}(i)|$ $%
\overset{p}{\rightarrow }$ $0$ in the proof of Theorem \ref{th:max_p} below.

\citet{Chernozhukov_etal2013,Chernozhukov_etal2016} work with normalized
sample means $1/\sqrt{n}\sum_{t=1}^{n}x_{t}(i)$, where $\{x_{t}(i)%
\}_{t=1}^{n}$ are independent, zero mean and square integrable random
variables. They bypass extreme value theoretic arguments and therefore do
not need to restrict dependence across $i$. This is accomplished by
developing new tools for deriving Gaussian approximations based on Slepian
and Sudakov-Fernique methods. They prove the impressive result that for some
$(K,\zeta )$ $>$ $0$, and $\mathcal{L}_{n}$ $\rightarrow $ $\infty $ with $%
\mathcal{L}_{n}$ $=$ $O(e^{o(n^{c})})$ for some $c$ $>$ $0$:
\begin{equation}
\mathcal{A}_{\mathcal{L}_{n},n}\equiv \sup_{c\geq 0}\left\vert P\left(
\max_{1\leq i\leq \mathcal{L}_{n}}\left\vert \frac{1}{\sqrt{n}}%
\sum_{t=1}^{n}x_{t}(i)\right\vert \leq c\right) -P\left( \max_{1\leq i\leq
\mathcal{L}_{n}}\left\vert \frac{1}{\sqrt{n}}\sum_{t=1}^{n}y_{t}(i)\right%
\vert \leq c\right) \right\vert \leq Kn^{-\zeta }  \label{approx_Guass}
\end{equation}%
where $y_{t}(i)$ are zero mean normally distributed with covariance function
$E[y_{t}(i)y_{t}(j)]$ $=$ \linebreak $1/n\sum_{t=1}^{n}E[x_{t}(i)x_{t}(j)]$.

The literature on such Gaussian couplings has a substantial history, where
generally $\mathcal{L}_{n}$ $=$ $O(n^{c})$ for some $c$. See, for example, %
\citet[Chapter 10]{Pollard2002} for a general review, see \cite%
{Yurinskii1977} for a seminal result, and see
\citet[Lemma
2.12]{DudleyPhilipp1983} and \cite{LeCam1988}. See also %
\citet{Portnoy1985,Portnoy1986}, \cite{Gotze1991} and
\citet[Appendix
L]{Chernozhukov_etal2013_sm} for theory and references on high dimensional
Gaussian central limit theory.

In the above literature, Gaussianicity is key, which therefore neglects
non-standard asymptotics, including heavy tailed data or non-stationary
data. It also does not apply in general when working with filtered data: in
this case an intermediate step is required linking two sequences in which
neither may be Gaussian. As such, it does not include cases where a filter
leads to non-standard asymptotics, including with some parameters are weakly
or non-identified and when a parameter boundary value occurs %
\citep[e.g.][]{Andrews1999,AndrewsCheng2012}.

Consider, for example, random functions $x_{t}(i,\theta _{0}(i))$ where $%
\theta _{0}(i)$ is not observed but estimable, write $\mathcal{X}%
_{n}(i,\theta _{0}(i))$ $\equiv $ 1/$\sqrt{n}\sum_{t=1}^{n}x_{t}(i,\theta
_{0}(i))$, and let $\hat{\theta}_{n}(i)$ be a plug-in estimator for $\theta
_{0}(i)$. Write $\mathcal{X}_{n}(i)$ $\equiv $ $\mathcal{X}_{n}(i,\hat{\theta%
}_{n}(i))$. In order to gain inference on $\max_{1\leq i\leq \mathcal{L}%
_{n}}|\mathcal{X}_{n}(i)|$ we first require an asymptotic expansion $%
\mathcal{X}_{n}(i)$ $=$ $\mathcal{Y}_{n}(i)$ $+$ $o_{p}(1)$ for some process
$\{\mathcal{Y}_{n}(i)\}$ that ostensibly depends on $\theta _{0}(i)$ and
pre-asymptotic properties of $\hat{\theta}_{n}(i)$, leading to $|\max_{1\leq
i\leq \mathcal{L}_{n}}|\mathcal{X}_{n}(i)|$ $-$ $\max_{1\leq i\leq \mathcal{L%
}_{n}}|\mathcal{Y}_{n}(i)||$ $\overset{p}{\rightarrow }$ $0$ as in (\ref%
{max|X-Y|}). Even if $\mathcal{Y}_{n}(i)$ is \textit{asymptotically} normal,
it need not be Gaussian, hence the Gaussian approximation literature does
not apply. Of course, $\mathcal{Y}_{n}(i)$ need not be asymptotically normal
for (\ref{max|X-Y|})\ to apply. In Section \ref{sec:ex} we illustrate
non-standard asymptotics by allowing for unit root non-stationarity, and we
treat a white noise test for an unobserved regression error term to
illustrate the use of a filter.

The maximum of an increasing sequence of normalized sample covariances $\hat{%
\gamma}_{n}(i)$ $\equiv $ $1/n\sum_{t=1}^{n}x_{t}x_{t-i}$ has been studied
at least since \cite{Berman1964} and \cite{Hannan1974}. See \cite{Jirak2011}
and \cite{XiaoWu2014} for recent theory and references. In this literature $%
x_{t}$ is assumed observed, the exact asymptotic distribution form of a
suitably normalized $\sqrt{n}\max_{1\leq i\leq \mathcal{L}_{n}}|\hat{\gamma}%
_{n}(i)$ $-$ $E[x_{t}x_{t-i}]|$ is sought, and underlying assumptions ensure
$\sqrt{n}(\hat{\gamma}_{n}(i)$ $-$ $E[x_{t}x_{t-i}])$ converges in finite
dimensional distributions to a Gaussian law $\mathcal{Z}(i)$. In turn, $%
\max_{1\leq i\leq \mathcal{L}_{n}}|\mathcal{Z}(i)|$ must converge in law to
a well defined random variable, which requires asymptotic independence $E[%
\mathcal{Z}(i)\mathcal{Z}(j)]$ $\rightarrow $ $0$ as $|i$ $-$ $j|$ $%
\rightarrow $ $\infty $. See, e.g., \citet[Chapter
6]{Leadbetter_etal1983}, \citet{Husler1986,Husler1993}, \cite%
{HombelMcCormick1995} and \citet[Chapter  9]{FalkHuslerReiss2011}. We
sidestep extreme value theoretic arguments, dependence across $i$ is
unrestricted, and residuals are allowed, hence $x_{t}$ need not be observed
in practice. See Section \ref{sec:ex}.

We cannot generally provide an upper bound on the divergence rate $\mathcal{L%
}_{n}$ $\rightarrow $ $\infty $, similar to ones in %
\citet{Chernozhukov_etal2013,Chernozhukov_etal2016} in the Gaussian coupling
literature, and \cite{XiaoWu2014} in the extreme value theory literature.
This is an unavoidable cost for our ($i$) basing probabilistic statements
like (\ref{max|X-Y|}) on a general array convergence result that itself does
not make use of probabilistic properties of (Gaussian) random variables; and
($ii $) allowing for filtered data and therefore requiring asymptotic
linkages between maxima that do not involve Gaussian processes. %
\citet{Chernozhukov_etal2013,Chernozhukov_etal2015,Chernozhukov_etal2016},
however, appear to have the sharpest result and most general bound on $%
\mathcal{L}_{n}$ $\rightarrow $ $\infty $\ for Gaussian approximations (for
independent data).

We do not treat a bootstrap theory, for example for the sample mean maximum
\linebreak $\max_{1\leq i\leq \mathcal{L}_{n}}|1/\sqrt{n}%
\sum_{t=1}^{n}x_{t}(i)|$, because by using our main results any existing
bootstrap theory will hold under its specified pointwise assumptions. See
\cite{HillMotegi2017}.\medskip

The remaining sections are organized as follows. In Section \ref{sec:main}
we present the main results (\ref{max|X|}) and (\ref{max|X-Y|}) of \cite%
{HillMotegi2017}. These are used in Section \ref{sec:filter} for convergence
of maxima when filtered residuals are used. Examples are provided in Section %
\ref{sec:ex}, and concluding remarks are left for Section \ref{sec:conclude}%
.\medskip

In the following $|\cdot |$ and $||\cdot ||$\ are respectively the $l_{1}$-
and $l_{2}$-matrix norms.

\section{Main Results\label{sec:main}}

All random variables are assumed to exist on a complete measure space, and
probabilities where applicable are outer probability measures. See, e.g., %
\citet[p. 101]{Dudley1984} and \citet[Appendix C]{Pollard1984}. See \cite%
{DudleyPhilipp1983} for theory that sidesteps measurability challenges
specifically for normalized means.

The first result due to \cite{HillMotegi2017}\ concerns convergence of
arrays. It is based on, and extends, a result from \citet[Lemma
1]{BoehmeRosenfeld1974}. They work in first countable spaces. Because we
require some structure on the space we work in, and because any metric space
is first countable, we simply work in $(\mathbb{R},d)$ with metric $d$ for
ease of notation. All proofs are placed in the appendix.
%\hyperref[app:proofs]{appendix}.

\begin{lemma}
\label{lm:array_conv}Assume the array $\{\mathcal{A}_{k,n}$ $:$ $1$ $\leq $ $%
k$ $\leq $ $\mathcal{I}_{n}\}_{n\geq 1}$ lies $(\mathbb{R},d)$, where $\{%
\mathcal{I}_{n}\}_{n\geq 1}$ is a sequence of positive integers, $\mathcal{I}%
_{n}$ $\rightarrow $ $\infty $ as $n$ $\rightarrow $ $\infty $. Let $%
\lim_{n\rightarrow \infty }\mathcal{A}_{k,n}$ $=$ $0$ for each fixed $k$,
and $\mathcal{A}_{k,n}$ $\leq $ $\mathcal{A}_{k+1,n}$ for each $n$ and all $%
k $. Then $\lim_{n\rightarrow \infty }\mathcal{A}_{\mathcal{L}_{n},n}$ $=$ $0
$ for some sequence $\{\mathcal{L}_{n}\}$ of positive integers, $\mathcal{L}%
_{n}$ $\rightarrow $ $\infty $ and $\mathcal{L}_{n}$ $\leq $ $\mathcal{I}%
_{n} $, that is not unique.
\end{lemma}

\begin{remark}
\normalfont\cite{BoehmeRosenfeld1974} only require $\lim_{k\rightarrow
\infty }\lim_{n\rightarrow \infty }\mathcal{A}_{k,n}$ $=$ $0$, and do not
impose monotonicity, thus $\lim_{n\rightarrow \infty }\mathcal{A}_{k,n}$ $=$
$0$ $\forall k$ and $\mathcal{A}_{k,n}$ $\leq $ $\mathcal{A}_{k+1,n}$ are
stronger assumptions. The first property, however, is key towards proving $%
\lim_{n\rightarrow \infty }\mathcal{A}_{\mathcal{L}_{n},n}$ $=$ $0$, rather
than merely a subsequence $\lim_{l\rightarrow \infty }\mathcal{A}_{\mathcal{L%
}_{n_{l}},n_{l}}$ $=$ $0$ as in \cite{BoehmeRosenfeld1974}. Monotonicity is
used to identify $\mathcal{L}_{n}$ as a function only of $n$ based on using
a multiple subsequence argument. The maximum over a subsequence of positive
values satisfies monotonicity, and $\lim_{n\rightarrow \infty }\mathcal{A}%
_{k,n}$ $=$ $0$ $\forall k$ holds when applied to pointwise probability
convergence problems discussed in the sequel.
\end{remark}

\begin{remark}
\normalfont Monotonicity implies $\{\mathcal{L}_{n}\}$ is not unique since $%
\lim_{n\rightarrow \infty }\mathcal{A}_{\mathcal{\mathring{L}}_{n},n}$ $=$ $%
0 $ for any sequence $\{\mathcal{\mathring{L}}_{n}\}$ of positive integers
with $\mathcal{\mathring{L}}_{n}$ $\rightarrow $ $\infty $ and $\lim
\sup_{n\rightarrow \infty }\{\mathcal{\mathring{L}}_{n}/\mathcal{L}_{n}\}$ $%
< $ $1$.
\end{remark}

The next result uses Lemma \ref{lm:array_conv} as the basis for deriving (%
\ref{max|X|}) and (\ref{max|X-Y|}).

\begin{theorem}
\label{th:max_p}Let $\{\mathcal{X}_{n}(i),\mathcal{Y}_{n}(i)$ $:$ $1$ $\leq $
$i$ $\leq $ $\mathcal{I}_{n}\}_{n\geq 1}$ be arrays of random variables,
where $\{\mathcal{I}_{n}\}_{n\geq 1}$ is a sequence of positive integers, $%
\mathcal{I}_{n}$ $\rightarrow $ $\infty $ as $n$ $\rightarrow $ $\infty $%
.\medskip \newline
$a.$ If $\mathcal{X}_{n}(i)$ $\overset{p}{\rightarrow }$ $0$ for each fixed $%
i$, then $\max_{1\leq i\leq \mathcal{L}_{n}}|\mathcal{X}_{n}(i)|$ $\overset{p%
}{\rightarrow }$ $0$ for some sequence $\{\mathcal{L}_{n}\}$ of positive
integers with $\mathcal{L}_{n}$ $\rightarrow $ $\infty $, that is not
unique.\medskip \newline
$b.$ If $\mathcal{X}_{n}(i)$ $-$ $\mathcal{Y}_{n}(i)$ $\overset{p}{%
\rightarrow }$ $0$ for each fixed $i$\ then for some sequence $\{\mathcal{L}%
_{n}\}$ of positive integers with $\mathcal{L}_{n}$ $\rightarrow $ $\infty $%
, that is not unique:
\begin{equation*}
\left\vert \max_{1\leq i\leq \mathcal{L}_{n}}\left\vert \mathcal{X}%
_{n}(i)\right\vert -\max_{1\leq i\leq \mathcal{L}_{n}}\left\vert \mathcal{Y}%
_{n}(i)\right\vert \right\vert \leq \max_{1\leq i\leq \mathcal{L}%
_{n}}\left\vert \mathcal{X}_{n}(i)-\mathcal{Y}_{n}(i)\right\vert \overset{p}{%
\rightarrow }0.
\end{equation*}
\end{theorem}

\begin{remark}
\normalfont A similar result exists under \emph{almost sure} convergence,
although a different argument is required. See the supplemental material
\citet[Appendix
B]{max_mean_supp_mat}.
\end{remark}

\begin{remark}
\normalfont The method of proof for Lemma \ref{lm:array_conv} shows the
\emph{existence}\ of such a sequence $\{\mathcal{L}_{n}\}$ and therefore
cannot provide an upper bound on the rate $\mathcal{L}_{n}$ $\rightarrow $ $%
\infty $. This seems unavoidable since we are working with general array
convergence rather than, for example, the specific attributes of Gaussian
probability tails. The payoff is that such generality ultimately permits
non-Gaussian couplings, a data filter, and arbitrary dependence across $i$,
as we treat in turn below.
\end{remark}

\section{Asymptotics for Maxima Based on Filtered Data\label{sec:filter}}

We now work with a parametric array $\{\mathcal{X}_{n}(i,\theta _{0}(i))$ $:$
$1$ $\leq $ $i$ $\leq $ $\mathcal{L}\}_{n\geq 1}$, where $\theta _{0}(i)$ is
an unknown but estimable parameter in $\mathbb{R}^{k}$ that may depend on $i$%
. Let $\hat{\theta}_{n}(i)$ be an estimator of $\theta _{0}(i)$. We assume $%
\mathcal{X}_{n}(i,\theta _{0}(i))$ is unobserved while $\mathcal{X}_{n}(i,%
\hat{\theta}_{n}(i))$ is observed. Our leading example is a sample serial
correlation coefficient for time series regression model errors (see Section %
\ref{sec:ex}). Another is a sample mean of an observed time series scaled by
its conditional variance, e.g. GARCH residuals.

Our primary goal is to prove under fairly general conditions that for some
stochastic process $\{\mathcal{Z}_{n}(i)\}_{i\in \mathbb{N}}$, and some
non-unique $\{\mathcal{L}_{n}\}$, $\mathcal{L}_{n}$ $\rightarrow $ $\infty $:%
\begin{equation}
\left\vert \max_{1\leq i\leq \mathcal{L}_{n}}\left\vert \mathcal{X}_{n}(i,%
\hat{\theta}_{n}(i))\right\vert -\max_{1\leq i\leq \mathcal{L}%
_{n}}\left\vert \mathcal{Z}_{n}(i)\right\vert \right\vert \overset{p}{%
\rightarrow }0.  \label{Xn_approx_max}
\end{equation}%
In view of Theorem \ref{th:max_p}, it suffices to prove $\mathcal{X}_{n}(i,%
\hat{\theta}_{n}(i))$ $-$ $\mathcal{Z}_{n}(i))$ $\overset{p}{\rightarrow }$ $%
0$ for each fixed $i$. Sufficient conditions follow. These are not the most
general possible, but give a reasonably general environment to work in.

\begin{assumption}
\label{assum:X_filter} $\ \ \ \medskip $\newline
$a$. $\mathcal{X}_{n}(i,\theta )$ is continuous and differentiable in $%
\theta $ $\in $ $\Theta (i)$, where $\Theta (i)$ is a compact subset of $%
\mathbb{R}^{k}$, $k$ $\in $ $\mathbb{N}$.\medskip \newline
$b.$ $\theta _{0}(i)$ $\in $ $\Theta (i)$.\medskip \newline
$c.$ There exists a continuous non-stochastic function $\mathcal{D}(i,\cdot
) $ $:$ $\Theta (i)$ $\rightarrow $ $\mathbb{R}^{k}$, and a compact
neighborhood $\mathcal{N}_{0}(i)$ of $\theta _{0}(i)$ with positive Lebesgue
measure such that%
\begin{equation*}
\sup_{\theta \in \mathcal{N}_{0}(i)}\left\vert \frac{1}{\sqrt{n}}\frac{%
\partial }{\partial \theta }\mathcal{X}_{n}(i,\theta )-\mathcal{D}(i,\theta
)\right\vert \overset{p}{\rightarrow }0.
\end{equation*}%
$d.$ There exists stochastic processes $\{\mathcal{S}_{n}(i),\mathcal{M}%
_{n}(i)$ $:$ $i$ $\in $ $\mathbb{N}\}$ such that $\mathcal{X}_{n}(i,\theta
_{0}(i))$ $=$ $\mathcal{S}_{n}(i)+$ $o_{p}(1)$ and $\sqrt{n}(\hat{\theta}%
_{n}(i)$ $-$ $\theta _{0}(i))$ $=$ $\mathcal{M}_{n}(i)$ $+$ $o_{p}(1)$ for
each $i$. Moreover, there exists non-degenerate stochastic processes $\{%
\mathcal{S}(i),\mathcal{M}(i)$ $:$ $i$ $\in $ $\mathbb{N}\}$ such that $(%
\mathcal{S}_{n}(i),\mathcal{M}_{n}(i))$ $\overset{d}{\rightarrow }$ $(%
\mathcal{S}(i),\mathcal{M}(i))$ for each $i$.
\end{assumption}

\begin{remark}
\normalfont Uniform convergence (c) holds when $(\partial /\partial \theta )%
\mathcal{X}_{n}(i,\theta )/\sqrt{n}$ $\overset{p}{\rightarrow }$ $\mathcal{D}%
(i,\theta )$ for each $\theta $ $\in $ $\mathcal{N}_{0}(i)$, and a
stochastic equicontinuity. condition holds. The latter holds, for example,
when $\mathcal{X}_{n}(i,\theta )$ is twice continuously differentiable and
the envelope%
\begin{equation*}
E\left[ \frac{1}{\sqrt{n}}\sup_{\theta \in \mathcal{N}_{0}(i)}\left\vert
\frac{\partial ^{2}}{\partial \theta _{i}\partial \theta _{j}}\mathcal{X}%
_{n}(i,\theta )\right\vert \right] =O(1).
\end{equation*}%
See, for example, \cite{Newey1991} and \cite{Andrews1992}.
\end{remark}

\begin{remark}
\normalfont Joint convergence (d) holds, for example, when $\mathcal{X}%
_{n}(i,\cdot )$ $=$ $1/\sqrt{n}\sum_{t=1}^{n}x_{t}(i)$ where $x_{t}(i)$ has
a zero mean and is square integrable, $\sqrt{n}(\hat{\theta}_{n}(i)$ $-$ $%
\theta _{0}(i))$ $=$ $1/\sqrt{n}\sum_{t=1}^{n}m_{t}(i)$ $+$ $o_{p}(1)$ for
some zero mean square integrable random variables $m_{t}(i)$, and $%
\{x_{t}(i),m_{t}(i)\}$ satisfy suitable moment and dependence properties. An
example concerning a residuals white noise test is provided in Section \ref%
{sec:ex}.
\end{remark}

\begin{remark}
\normalfont Conditions (b) and (d) allow for non-standard cases. One example
is when $\theta _{0}(i)$ lies on the boundary of $\Theta (i)$ %
\citep[e.g.][]{Andrews1999}, and another is when a component of $\theta
_{0}(i)$\ is weakly or non-identified \citep[e.g.][]{AndrewsCheng2012}. In
each case $\sqrt{n}(\hat{\theta}_{n}(i)$ $-$ $\theta _{0}(i))$ $=$ $\mathcal{%
M}_{n}(i)$ $+$ $o_{p}(1)$ $\overset{d}{\rightarrow }$ $\mathcal{M}(i))$ is
non-Gaussian.
\end{remark}

\begin{theorem}
\label{th:filter_X}Let Assumption \ref{assum:X_filter} hold, and write $%
\mathcal{D}(i)$ $\equiv $ $\mathcal{D}(i,\theta _{0})$ where $\mathcal{D}%
(i,\theta )$ is defined in Assumption \ref{assum:X_filter}.c. For some
sequence $\{\mathcal{L}_{n}\}$ of positive integers that is not unique, with
$\mathcal{L}_{n}$ $\rightarrow $ $\infty $: $|\max_{1\leq i\leq \mathcal{L}%
_{n}}|\mathcal{X}_{n}(i,\hat{\theta}_{n}(i))|$ $-$ $\max_{1\leq i\leq
\mathcal{L}_{n}}|\mathcal{S}_{n}(i)$ $+$ $\mathcal{D}(i)^{\prime }\mathcal{M}%
_{n}(i)||$ $\overset{p}{\rightarrow }$ $0$ and $|\max_{1\leq i\leq \mathcal{L%
}_{n}}|\mathcal{S}_{n}(i)+\mathcal{D}(i)^{\prime }\mathcal{M}_{n}(i)|$ $-$ $%
\max_{1\leq i\leq \mathcal{L}_{n}}|\mathcal{S}(i)+\mathcal{D}(i)^{\prime }%
\mathcal{M}(i)||$ $\overset{p}{\rightarrow }$ $0$, hence $|\max_{1\leq i\leq
\mathcal{L}_{n}}|\mathcal{X}_{n}(i,\hat{\theta}_{n}(i))|$ $-$ $\max_{1\leq
i\leq \mathcal{L}_{n}}|\mathcal{S}(i)$ $+$ $\mathcal{D}(i)^{\prime }\mathcal{%
M}(i)||$ $\overset{p}{\rightarrow }$ $0.$\qquad
\end{theorem}

\section{Illustrations\label{sec:ex}}

We now consider maximum statistics that involve unit root test statistics
and a residuals white noise test statistic. Both require a non-Gaussian
approximation theory, demonstrating the unique applicability of Theorem \ref%
{th:max_p}.

\subsection{Unit Root Tests}

Consider unit root tests over a set of processes $\{y_{t}(i)$ $:$ $i$ $\in $
$\mathbb{N}\}$. Suppose $y_{t}(i)$ $=$ $\phi _{0}(i)y_{t-1}(i)$ $+$ $%
\epsilon _{t}(i)$, $|\phi _{0}(i)|$ $\leq $ $1$, for $i$ $\in $ $\mathbb{N}$%
. We assume $\{\epsilon _{t}(i)$ $:$ $i$ $\in $ $\mathbb{N}\}$ lies in a
probability measure space $(\Omega ,\sigma (\cup _{t\in \mathbb{N}}\mathcal{F%
}_{t}),\mathcal{P})$, where $\{\mathcal{F}_{t}\}_{t\in \mathbb{N}}$ is a
sequence of $\sigma $-felds. Assume $\epsilon _{t}(i)$ is $\mathcal{F}_{t}$%
-measurable for each $i$. Write $\mathcal{F}_{s}^{t}$ $\equiv $ $\sigma
(\cup _{\tau =s}^{t}\mathcal{F}_{\tau })$. Zero mean $\epsilon _{t}(i)$ is
stationary, and $E|\epsilon _{t}(i)|^{r}$ $<$ $\infty $ for some $r$ $>$ $2$%
. Moreover, $\alpha $-mixing coefficients $\alpha _{h}$ $\equiv $ $\sup_{%
\mathcal{A}\in \mathcal{F}_{-\infty }^{t-h},\mathcal{B}\in \mathcal{F}%
_{t}^{\infty }}|\mathcal{P}(\mathcal{A}\cap \mathcal{B})$ $-$ $\mathcal{P}(%
\mathcal{A})P(\mathcal{B})|$ satisfy $\alpha _{h}$ $=$ $O(h^{r/(r-2)}/\ln
(h))$. Thus, $\epsilon _{t}(i)$ and any measurable function of $\{\epsilon
_{t}(i)$ $:$ $1$ $\leq $ $i$ $\leq $ $k\}$ has coefficients $\alpha _{h}$.
Let $\hat{\phi}_{n}(i)$ be the least squares estimator.

We want to test the hypothesis that all processes have a unit root $H_{0}$ $%
: $ $\phi _{0}(i)$ $=$ $1$ for all $i$ $\in $ $\mathbb{N}$. The proposed
test statistic is
\begin{equation*}
\mathcal{T}_{n}(k)\equiv n\max_{1\leq i\leq k}\left\vert \hat{\phi}%
_{n}(i)-1\right\vert .
\end{equation*}

Define $\sigma ^{2}(i)$ $\equiv $ $\lim_{n\rightarrow \infty }E[(1/\sqrt{n}%
\sum_{t=1}^{n}\epsilon _{t}(i))^{2}]$ and $\sigma _{\epsilon }^{2}(i)$ $%
\equiv $ $E[\epsilon _{t}^{2}(i)]$, and define
\begin{equation*}
\mathcal{T}(i)\equiv \frac{1}{2}\left\{ \mathcal{W}(i,1)^{2}-\sigma
_{\epsilon }^{2}(i)/\sigma ^{2}(i)\right\} /\int_{0}^{1}\mathcal{W}(i,\nu
)^{2}d\nu
\end{equation*}%
where $\{\mathcal{W}(\cdot ,\nu )$ $:$ $[0,1]\}$ are standard Wiener
processes. $\mathcal{T}(i)$ is the well known limit law for the least
squares estimator when there is a unit root and $\epsilon _{t}(i)$ is
possibly dependent: $n(\hat{\phi}_{n}(i)$ $-$ $1)$ $\overset{d}{\rightarrow }
$ $\mathcal{T}(i)$ \citep[Theorem 3.1]{Phillips1987}. If $\epsilon _{t}(i)$
is iid then $\sigma _{\epsilon }^{2}(i)/\sigma ^{2}(i)$ $=$ $1$, cf. \cite%
{White1958}.

Apply the mapping theorem to yield $n|\hat{\phi}_{n}(i)$ $-$ $1|$ $\overset{d%
}{\rightarrow }$ $|\mathcal{T}(i)|$ for each $i$ under $H_{0}$. Therefore $n|%
\hat{\phi}_{n}(i)$ $-$ $1|$ $-$ $|\mathcal{T}(i)|$ $=$ $o_{p}(1)$ for each $%
i $ under $H_{0}$. Now invoke Theorem \ref{th:max_p} to deduce $|\mathcal{T}%
_{n}(\mathcal{L}_{n})$ $-$ $\max_{1\leq i\leq \mathcal{L}_{n}}|\mathcal{T}%
(i)||$ $\overset{p}{\rightarrow }$ $0$ under $H_{0}$\ for some non-unique
sequence of positive integers $\{\mathcal{L}_{n}\}$, $\mathcal{L}_{n}$ $%
\rightarrow $ $\infty $. Conversely, if $|\phi _{0}(i^{\ast })|$ $<$ $1$ for
some $i^{\ast }$ $\in $ $\mathbb{N}$, then such $n(\phi _{0}(i^{\ast })$ $-$
$1)$ $\rightarrow $ $-\infty $ will eventually dominate $\mathcal{T}%
_{n}(k)=\max_{1\leq i\leq k}|n(\hat{\phi}_{n}(i)$ $-$ $\phi _{0}(i))$ $-$ $%
n(\phi _{0}(i)$ $-$ $1)|$ for any $k$ $\geq $ $i^{\ast }$. Therefore $%
\mathcal{T}_{n}(\mathcal{L}_{n})$ $\overset{p}{\rightarrow }$ $\infty $ for
any sequence of positive integers $\{\mathcal{L}_{n}\}$, $\mathcal{L}_{n}$ $%
\rightarrow $ $\infty $.

The limit law $\mathcal{T}(i)$ is not pivotal since it contains nuisance
parameters $\sigma _{\epsilon }^{2}(i)$ and $\sigma ^{2}(i)$. Let $\hat{%
\sigma}_{n,\epsilon }^{2}(i)$ and $\hat{\sigma}_{n}^{2}(i)$ be consistent
estimators of $\sigma _{\epsilon }^{2}(i)$ and $\sigma ^{2}(i)$ respectively
\citep[see][Section
4]{Phillips1987}, and define%
\begin{equation*}
\mathcal{\tilde{T}}_{n}(k)\equiv \max_{1\leq i\leq k}\left\vert n\left( \hat{%
\phi}_{n}(i)-\frac{\frac{1}{n}\frac{1}{2}\left( \hat{\sigma}_{n}^{2}(i)-\hat{%
\sigma}_{n,\epsilon }^{2}(i)\right) }{1/n^{2}\sum_{i=2}^{n}y_{t-1}^{2}(i)}%
-1\right) \right\vert =\max_{1\leq i\leq k}\left\vert \mathcal{\tilde{T}}%
_{n,i}\right\vert ,
\end{equation*}%
where $\mathcal{\tilde{T}}_{n,i}$\ is implicitly defined. The term $\mathcal{%
\tilde{T}}_{n,i}$ was proposed in{\ \citet[Section 5]{Phillips1987} as an
adjustment that leads to a pivotal asymptotic law.}

Denote $\mathcal{\tilde{T}}(i)$ $\equiv $ $(1/2)(\mathcal{W}(i,1)^{2}$ $-$ $%
1)/\int_{0}^{1}\mathcal{W}(i,\nu )^{2}d\nu $, the now classic limit law for
the least squares estimator when there is a unit root and iid error %
\citep{White1958}. Theorem 5.1 in \cite{Phillips1987} implies $\mathcal{%
\tilde{T}}_{n,i}$ $=$ $\mathcal{\tilde{T}}(i)$ $+$ $o_{p}(1)$. Now use
Theorem \ref{th:max_p} to yield $|\mathcal{\tilde{T}}_{n}(\mathcal{L}_{n})$ $%
-$ $\max_{i\in \mathcal{I}(\mathcal{L}_{n})}\mathcal{\tilde{T}}(i)|$ $%
\overset{p}{\rightarrow }$ $0$ under $H_{0}$ for some non-unique sequence of
positive integers $\{\mathcal{L}_{n}\}$, $\mathcal{L}_{n}$ $\rightarrow $ $%
\infty $.

\subsection{Residuals White Noise Test}

Consider a process $\{y_{t}\}$ modeled as an AR($p$) for finite $p$ $\geq $ $%
1$,
\begin{equation*}
y_{t}=c_{0}+\sum_{i=1}^{p}\phi _{0,i}y_{t-1}+\epsilon _{t}=\theta
_{0}^{\prime }x_{t}+\epsilon _{t}
\end{equation*}%
where $\theta _{0}$ $\equiv $ $[c_{0},\phi _{0}^{\prime }]^{\prime }$, $%
x_{t} $ $\equiv $ $[1,y_{t-1},...,y_{t-p}]^{\prime }$, and $1$ $-$ $%
\sum_{i=1}^{p}\phi _{0,i}z^{i}$ has roots outside the unit circle. The AR
model is assumed to be pseudo true in the sense that $\theta _{0}$ is the
unique point in the interior of compact $\Theta $ that satisfies $E[\epsilon
_{t}x_{t}]$ $=$ $0$. Assume $E[y_{t}^{2}]$ $>$ $0$ and $E|y_{t}|^{r}$ $<$ $%
\infty $ for some $r$ $>$ $4$. Define $\sigma $-fields $\mathcal{F}_{t}$ $%
\equiv $ $\sigma (y_{\tau }$ $:$ $\tau $ $\leq $ $t)$ and $\mathcal{F}%
_{s}^{t}$ $\equiv $ $\sigma (y_{\tau }$ $:$ $s$ $\leq $ $\tau $ $\leq $ $t)$%
, and assume $\mathcal{F}_{t-1}$ $\subset $ $\mathcal{F}_{t}$ $\forall t$.
Assume $y_{t}$ is stationary $\alpha $-mixing with coefficients $\alpha _{h}$
$\equiv $ $\sup_{\mathcal{A}\subset \mathcal{F}_{-\infty }^{t-m},\mathcal{B}%
\subset \mathcal{F}_{t}^{\infty }}|P(\mathcal{A}\cap \mathcal{B})$ $-$ $P(%
\mathcal{A})P(\mathcal{B})|$ $=$ $O(h^{r/(r-2)}/\ln (h))$. Sufficient
conditions for the strong mixing property in linear processes are presented
in \cite{Gorodetskii1977} and \cite{Withers1981}, amongst others.

We want to test $H_{0}$ $:$ $E[\epsilon _{t}\epsilon _{t-h}]$ $=$ $0$ $%
\forall h$ $\geq $ $1$. Let $\hat{\theta}_{n}$ $\equiv $ $[\hat{c}_{n},\hat{%
\phi}_{n}^{\prime }]$ be the least squares estimator of $\theta _{0}$, and
define an error function, and the residual sample serial correlation at lag $%
h$ $\geq $ $1$:
\begin{equation*}
\epsilon _{t}(\theta )\equiv y_{t}-\theta ^{\prime }x_{t}\text{ \ and \ }%
\mathcal{X}_{n}(i,\theta )\equiv \sqrt{n}\frac{1/n\sum_{t=1+h}^{n}\epsilon
_{t}(\theta )\epsilon _{t-h}(\theta )}{1/n\sum_{t=1}^{n}\epsilon
_{t}^{2}(\theta )}.
\end{equation*}%
A valid test can be based on the maximum absolute correlation, $\max_{1\leq
h\leq \mathcal{L}_{n}}|\mathcal{X}_{n}(i,\hat{\theta}_{n})|$, cf. \cite%
{HillMotegi2017}. We only provide a proof linking $\max_{1\leq h\leq
\mathcal{L}_{n}}|\mathcal{X}_{n}(i,\hat{\theta}_{n})|$ to the maximum of a
process that depends on $\epsilon _{t}$ and properties of the plug-in
estimator, based on an asymptotic expansion. See
\citet[Section 2, Theorem
2.5]{HillMotegi2017} for an asymptotically valid dependent wild bootstrap
based on the expansion.

The following is based on arguments in \citet[Lemma 2.1]{HillMotegi2017}.%
\footnote{%
As in the present example, \cite{HillMotegi2017} work with a residuals
sample correlation in a general parametric regression model setting. Since
the form of $\mathcal{X}_{n}(i,\theta )$\ is therefore known, they do not
need to verify conditions like Assumption \ref{assum:X_filter}, and instead
deliver an expansion using a more direct proof.} The proof relies on a
standard expansion, and Theorem \ref{th:max_p}. We therefore prove the claim
in the supplemental material \cite{max_mean_supp_mat}. Define%
\begin{eqnarray*}
&&\mathcal{D}(h)=-E\left[ \epsilon _{t}x_{t-h}\right] -E\left[ \epsilon
_{t-h}x_{t}\right] \text{ and }z_{t}(h)\equiv \frac{1}{E\left[ \epsilon
_{t}^{2}\right] }\epsilon _{t}\epsilon _{t-h}-\mathcal{D}(h)^{\prime }\left(
E\left[ x_{t}x_{t}^{\prime }\right] \right) ^{-1}x_{t}\epsilon _{t} \\
&&w_{t}(\lambda )=\lambda _{1}\frac{1}{E\left[ \epsilon _{t}^{2}\right] }%
\epsilon _{t-h}+\lambda _{2}\left( E\left[ x_{t}x_{t}^{\prime }\right]
\right) ^{-1}x_{t}\text{ for arbitrary }\left( \lambda _{1},\lambda
_{2}\right) \in \mathbb{R}\text{, }\lambda _{1}^{2}+\lambda _{2}^{2}=1.
\end{eqnarray*}

\begin{theorem}
\label{th:wn} \ \ \ \medskip \newline
$a$. Let $H_{0}$ hold and assume $\lim_{n\rightarrow \infty }E(1/\sqrt{n}%
\sum_{t=1}^{n}w_{t}(\lambda )\epsilon _{t})^{2}]$ $>$ $0$. Moreover, assume
for any asymptotic draw $\{y_{t},\epsilon _{t}\}_{t=1}^{\infty }$ that $%
\inf_{\theta \in \Theta }|\epsilon _{t}(\theta )|$ $\geq $ $\iota $ $a.s.$
for all $t$ $\in $ $\mathbb{N}/S$ where $S$ is finite, and for some
non-random $\iota $ $>$ $0$. Assumption \ref{assum:X_filter} applies, and
therefore the conclusions of Theorem 3.1 hold. In particular, for some
non-unique sequence of positive integers $\{\mathcal{L}_{n}\}$, $\mathcal{L}%
_{n}$ $\rightarrow $ $\infty $:
\begin{equation}
\left\vert \max_{1\leq h\leq \mathcal{L}_{n}}\left\vert \mathcal{X}_{n}(i,%
\hat{\theta}_{n})\right\vert -\max_{1\leq h\leq \mathcal{L}_{n}}\left\vert
\frac{1}{\sqrt{n}}\sum_{t=1+h}^{n}z_{t}(h)\right\vert \right\vert \overset{p}%
{\rightarrow }0.  \label{Xz}
\end{equation}%
Moreover, $\max_{1\leq h\leq \mathcal{L}_{n}}||1/\sqrt{n}%
\sum_{t=1+h}^{n}z_{t}(h)|$ $-$ $\max_{1\leq h\leq \mathcal{L}_{n}}|\mathcal{Z%
}(h)||$ $\overset{p}{\rightarrow }$ $0$ where $\{\mathcal{Z}(h)$ $:$ $h$ $%
\in $ $\mathbb{N}\}$ is a zero mean Gaussian process with covariance kernel $%
E[\mathcal{Z}(h)\mathcal{Z}(\tilde{h})]$ $=$ $\lim_{n\rightarrow \infty
}1/n\sum_{s,t=1}E[z_{s}(h)z_{t}(\tilde{h})]$ and variance $E[\mathcal{Z}%
(h)^{2}]$ $<$ $\infty $.$\medskip $\newline
$b$. If $H_{0}$ is false then $\max_{1\leq h\leq \mathcal{L}_{n}}|\mathcal{X}%
_{n}(i,\hat{\theta}_{n})|$ $\overset{p}{\rightarrow }$ $\infty $ for any
sequence of positive integers $\{\mathcal{L}_{n}\}$, $\mathcal{L}_{n}$ $%
\rightarrow $ $\infty $.
\end{theorem}

\begin{remark}
\normalfont The bound $\lim_{n\rightarrow \infty }E(1/\sqrt{n}%
\sum_{t=1}^{n}w_{t}(\lambda )\epsilon _{t})^{2}]$ $>$ $0$ ensures a
non-degenerate limit theory for a key joint process arising in a sample
correlation first order expansion. The limit is finite by the mixing
property and $E|w_{t}(\lambda )\epsilon _{t}|^{r/2}$ $<$ $\infty $ where $r$
$>$ $4$\ by assumed $L_{r}$-boundedness \citep[Theorem 1.7]{Ibrag1975}. The
assumption $\inf_{\theta \in \Theta }|\epsilon _{t}(\theta )|$ $\geq $ $%
\iota $ $>$ $0$ $a.s.$ for all $t$ $\in $ $\mathbb{N}/S$ and finite $S$
expedites the expansion proof. This is mild since $\mathcal{F}_{t-1}$ $%
\subset $ $\mathcal{F}_{t}$ $\forall t$ implies $\inf_{\theta \in \Theta
}|\epsilon _{t}(\theta )|$ $>$ $0$ $a.s.$ $\forall t$ (see the proof of
Theorem \ref{th:wn}).
\end{remark}

\begin{remark}
\normalfont The result significantly augments known results in the
max-correlation literature by permitting residuals, and without restricting
dependence in the limit process $\{\mathcal{Z}(h)\}$. Existing extreme value
theory works with observed data, and imposes conditions that ensure $\sqrt{n}%
(\hat{\gamma}_{n}(i)$ $-$ $E[x_{t}x_{t-i}])$ converges in finite dimensional
distributions to a Gaussian law $\mathcal{Z}(i)$, and $\max_{1\leq i\leq
\mathcal{L}_{n}}|\mathcal{Z}(i)|$ converges in law to a well defined random
variable \citep[e.g.][]{XiaoWu2014}. The latter requires asymptotic
independence $E[\mathcal{Z}(i)\mathcal{Z}(j)]$ $\rightarrow $ $0$ as $|i$ $-$
$j|$ $\rightarrow $ $\infty $, cf. \citet[Chapter
6]{Leadbetter_etal1983} and \citet{Husler1986,Husler1993}. The high
dimensional central limit theory literature can tackle $\max_{1\leq h\leq
\mathcal{L}_{n}}|1/\sqrt{n}\sum_{t=1+h}^{n}z_{t}(h)|$, e.g. \cite%
{Chernozhukov_etal2013}, but not the intermediate step (\ref{Xz}) since $%
z_{t}(h)$ is generally not Gaussian.
\end{remark}

\section{Conclusion\label{sec:conclude}}

We provide a general result for convergence of arrays that permits a new
theory for the maximum of an increasing sequence of random variables. When
linking the maximum of two random variables $\mathcal{X}_{n}(i)$ and $%
\mathcal{Y}_{n}(i)$, unlike the extreme value theory and high dimensional
Gaussian central limit theory literatures, we do not require normality or
even asymptotic normality of $\mathcal{Y}_{n}(i)$. This permits new results
for maxima, covering heavy tailed data, non-stationary data, and filtered
data where the filter may lead to non-standard asymptotics. Two
illustrations are provided covering unit root tests and a residual white
noise test, both of which appear to be new. A shortcoming of our general
approach, based ultimately on \cite{Ramsey1930} theory and its implications
for array convergence, is that we cannot bound the allowed array dimension $%
\mathcal{L}_{n}$ as $n$ $\rightarrow $ $\infty $. This runs contrary to the
max-correlation literature \citep[e.g.][]{XiaoWu2014}, and the Gaussian
coupling literature, most recently punctuated by %
\citet{Chernozhukov_etal2013,Chernozhukov_etal2016}, where the best known
bounds on $\mathcal{L}_{n}$ are available.

\setcounter{equation}{0} \renewcommand{\theequation}{{\thesection}.%
\arabic{equation}} \appendix

\section{Appendix: Proofs\label{app:proofs}}

\noindent \textbf{Proof of Lemma \ref{lm:array_conv}.}\qquad We prove in
Step 1 that $\lim_{l\rightarrow \infty }\mathcal{A}_{\mathcal{L}%
(n_{l}),n_{l}}$ $=$ $0$ for some sequence of positive integers $%
\{n_{l}\}_{l=1}^{\infty }$, $n_{l}$ $<$ $n_{l+1}$ $\forall l$, and some
mapping $\mathcal{L}(n_{l})$ $\leq $ $\mathcal{L}(n_{l+1})$, $\mathcal{L}%
(n_{l})$ $\rightarrow $ $\infty $ and $n_{l}$ $\rightarrow $ $\infty $ as $l$
$\rightarrow $ $\infty $. We use that result in Step 2 to prove the
claim.\medskip \newline
\textbf{Step 1.}\qquad We now prove $\lim_{l\rightarrow \infty }\mathcal{A}_{%
\mathcal{L}(n_{l}),n_{l}}$ $=$ $0$. By assumption $\{\mathcal{A}_{k,n}$ $:$ $%
1$ $\leq $ $k$ $\leq $ $\mathcal{I}_{n}\}_{n\geq 1}$ lies in $(\mathbb{R},d)$%
, which is a first countable topological space, and $\lim_{k\rightarrow
\infty }\lim_{n\rightarrow \infty }\mathcal{A}_{k,n}$ $=$ $0$. Therefore, by
Lemma 1 in \cite{BoehmeRosenfeld1974} there exists a sequence of positive
integers $\{\mathcal{L}_{i}\}_{i=1}^{\infty }$, $\mathcal{L}_{i}$ $%
\rightarrow $ $\infty $ as $i$ $\rightarrow $ $\infty $, and an integer
mapping $n(\mathcal{L})$ $\rightarrow $ $\infty $ as $\mathcal{L}$ $%
\rightarrow $ $\infty $ such that $\lim_{i\rightarrow \infty }\mathcal{A}_{%
\mathcal{L}_{i},n(\mathcal{L}_{i})}$ $=$ $0$. The relation $n(\mathcal{L})$ $%
\rightarrow $ $\infty $ as $\mathcal{L}$ $\rightarrow $ $\infty $ holds by
construction of the array $\{\mathcal{A}_{k,n}$ $:$ $1$ $\leq $ $k$ $\leq $ $%
\mathcal{I}_{n}\}_{n\geq 1}$ with $\mathcal{I}_{n}$ $\rightarrow $ $\infty $
as $n$ $\rightarrow $ $\infty $.

We can always assume monotonicity: $\mathcal{L}_{i}$ $\leq $ $\mathcal{L}%
_{i+1}$ $\forall i$. Simply note that $\lim_{i\rightarrow \infty }\mathcal{A}%
_{\mathcal{L}_{i},n(\mathcal{L}_{i})}$ $=$ $0$ implies $\lim_{l\rightarrow
\infty }\mathcal{A}_{\mathcal{L}_{i_{l}},n(\mathcal{L}_{i_{l}})}$ $=$ $0$
for every infinite subsequence $\{i_{l}\}_{l\geq 1}$ of $\{i\}_{i\geq 1}$.
Since $\mathcal{L}_{i}\rightarrow \infty $ as $i$ $\rightarrow $ $\infty $,
we can find a subsequence $\{i_{l}^{\ast }\}_{l\geq 1}$ such that $%
i_{l}^{\ast }$ $\leq $ $i_{l+1}^{\ast }$ and $\mathcal{L}_{i_{l}^{\ast
}}\leq \mathcal{L}_{i_{l+1}^{\ast }}$ for each $l$. This follows from the
monotone subsequence theorem, which itself follows from Ramsey's (%
\citeyear{Ramsey1930}) theorem, cf. \cite{ErdosSzekeres1935} and \cite%
{BurkillMirsky1973}. Now define $\mathcal{L}_{l}^{\ast }$ $\equiv $ $%
\mathcal{L}_{i_{l}^{\ast }}$, hence $\lim_{l\rightarrow \infty }\mathcal{A}_{%
\mathcal{L}_{l}^{\ast },n(\mathcal{L}_{l}^{\ast })}$ $=$ $0$ where $\mathcal{%
L}_{l}^{\ast }$ $\leq $ $\mathcal{L}_{l+1}^{\ast }$ and $\mathcal{L}%
_{l}^{\ast }$ $\rightarrow $ $\infty $ as $l$ $\rightarrow $ $\infty $.

Now let $\{n_{i}\}_{i=1}^{\infty }$ and $\{\mathcal{L}(n_{i})\}_{i=1}^{%
\infty }$\ be any sequences satisfying $n_{i}$ $=$ $n(\mathcal{L}_{i})$ and $%
\mathcal{L}(n_{i})$ $=$ $\mathcal{L}_{i}$. Hence by the above argument $%
\mathcal{L}(n_{i})$ $\leq $ $\mathcal{L}(n_{i}$ $+$ $1)$, $\mathcal{L}%
(n_{i}) $ $\rightarrow $ $\infty $ and $n_{i}$ $\rightarrow $ $\infty $,
such that $\lim_{i\rightarrow \infty }\mathcal{A}_{\mathcal{L}(n_{i}),n_{i}}$
$=$ $0.$ Note that $\lim_{i\rightarrow \infty }\mathcal{A}_{\mathcal{L}%
(n_{i}),n_{i}}$ $=$ $0$ \textit{if and only if} $\lim_{l\rightarrow \infty }%
\mathcal{A}_{\mathcal{L}(n_{i_{l}}),n_{i_{l}}}$ $=$ $0$\ for every
subsequence $\{n_{i_{l}}\}_{l=1}^{\infty }$ of $\{n_{i}\}_{i=1}^{\infty }$.
Since $n_{i}$ $\rightarrow $ $\infty $ as $i$ $\rightarrow $ $\infty $, by
the monotone subsequence theorem there exists a strictly monotonically
increasing subsequence $\{n_{i_{l}}\}_{l=1}^{\infty }$. Therefore, as
required $\lim_{l\rightarrow \infty }\mathcal{A}_{\mathcal{L}(n_{l}),n_{l}}$
$=$ $0$ for some sequence of positive integers $\{n_{l}\}_{l=1}^{\infty }$, $%
n_{l}$ $<$ $n_{l+1}$ $\forall l$, and $\mathcal{L}(n_{l})$ $\leq $ $\mathcal{%
L}(n_{l+1})$, $\mathcal{L}(n_{l})$ $\rightarrow $ $\infty $ and $n_{l}$ $%
\rightarrow $ $\infty $ as $l$ $\rightarrow $ $\infty $.\medskip \newline
\textbf{Step 2.}\qquad By assumption $\lim_{n\rightarrow \infty }\mathcal{A}%
_{k,n}$ $=$ $0$ $\forall k$. Therefore:
\begin{equation}
\lim_{s\rightarrow \infty }\mathcal{A}_{k,n_{s}}=0\text{ for every }k\text{
and every infinite subsequence }\left\{ n_{s}\right\} _{s\geq 1}.
\label{limAkns}
\end{equation}%
Now repeat the Step 1 argument for each $\{\mathcal{A}_{k,n_{s}}\}_{s\geq 1}$%
: there exists a strictly monotonically increasing subsequence of positive
integers $\left\{ n_{s_{l}}\right\} _{l\geq 1}$ and some integer mapping $%
\mathcal{L}_{s}(n_{s_{l}})$ that may depend on $s$, with $n_{s_{l}}$ $%
\rightarrow $ $\infty $ and $\mathcal{L}_{s}(n_{s_{l}})$ $\rightarrow $ $%
\infty $\ as $l$ $\rightarrow $ $\infty $ $\forall s$, such that $%
\lim_{l\rightarrow \infty }\mathcal{A}_{\mathcal{L}%
_{s}(n_{s_{l}}),n_{s_{l}}} $ $=$ $0$ $\forall s$. As above, we may take $%
\mathcal{L}_{s}(\cdot )$ to be monotonic: $\mathcal{L}_{s}(\tilde{n})$ $\leq
$ $\mathcal{L}_{s}(\tilde{n}$ $+$ $1)$ $\forall \tilde{n}$.

Since monotonic $\mathcal{L}_{s}(\tilde{n})$ $\rightarrow $ $\infty $\ as $%
\tilde{n}$ $\rightarrow $ $\infty $ $\forall s$, there exists an integer
mapping $\mathcal{L}(\cdot )$ such that $\mathcal{L}(n)$ $\rightarrow $ $%
\infty $ as $n$ $\rightarrow $ $\infty $\ and for each $s$, $%
\limsup_{n\rightarrow \infty }\{\mathcal{L}(n)/\mathcal{L}_{s}(n)\}$ $<$ $1.$
By monotonicity $\mathcal{A}_{k,n}$ $\leq $ $\mathcal{A}_{k+1,n}$ this
mapping satisfies%
\begin{equation}
\lim_{l\rightarrow \infty }\mathcal{A}_{\mathcal{L}(n_{s_{l}}),n_{s_{l}}}%
\leq \lim_{l\rightarrow \infty }\mathcal{A}_{\mathcal{L}%
_{s}(n_{s_{l}}),n_{s_{l}}}=0\text{ }\forall s.  \label{limAlimA}
\end{equation}%
Notice $\mathcal{L}(\cdot )$ is not unique: for any $\mathcal{L}(\cdot )$
that satisfies (\ref{limAlimA}) there exists $\mathcal{\tilde{L}}(n)$ $%
\rightarrow $ $\infty $ such that \linebreak $\lim \sup_{n\rightarrow \infty
}\mathcal{\tilde{L}}(n)/\mathcal{L}(n)$ $<$ $1$, hence by monotonicity $%
\lim_{l\rightarrow \infty }\mathcal{A}_{\mathcal{\tilde{L}}%
(n_{s_{l}}),n_{s_{l}}}$ $\leq $ $\lim_{l\rightarrow \infty }\mathcal{A}_{%
\mathcal{L}(n_{s_{l}}),n_{s_{l}}}$ $=$ $0$.

Now write $\mathcal{B}_{n}$ $\equiv $ $\mathcal{A}_{\mathcal{L}(n),n}.$ By a
direct implication of (\ref{limAkns}) and (\ref{limAlimA}), for every
subsequence $\{\mathcal{B}_{n_{s}}\}_{s\geq 1}$ there exists a further
subsequence $\{\mathcal{B}_{n_{s_{l}}}\}_{l\geq 1}$ that converges $%
\lim_{l\rightarrow \infty }\mathcal{B}_{n_{s_{l}}}$ $=$ $0$. Therefore $%
\lim_{n\rightarrow \infty }\mathcal{B}_{n}$ $=$ $0$
\citep[see][p.
39]{Royden1988}. This proves $\lim_{n\rightarrow \infty }\mathcal{A}_{%
\mathcal{L}_{n},n}$ $=$ $0$ with $\mathcal{L}_{n}$ $=$ $\mathcal{L}(n)$ as
required. $\mathcal{QED}$.\bigskip \newline
\textbf{Proof of Theorem \ref{th:max_p}}.\medskip \newline
\textbf{Claim (a).}\qquad By assumption each $\mathcal{X}_{n}(i)\overset{p}{%
\rightarrow }0$ and therefore $\max_{1\leq i\leq k}|\mathcal{X}_{n}(i)|$ $%
\overset{p}{\rightarrow }$ $0$. Define $\mathcal{A}_{k,n}$ $\equiv $ $1$ $-$
$\exp \{-\max_{1\leq i\leq k}|\mathcal{X}_{n}(i)|)\}$ and $\mathcal{P}_{k,n}$
$\equiv $ $\int_{0}^{\infty }P(\mathcal{A}_{k,n}$ $>$ $\epsilon )d\epsilon $%
. By construction $\mathcal{A}_{k,n}$ $\in $ $[0,1]$ $a.s$. $\forall k$.
Lebesgue's dominated convergence theorem, and $\mathcal{A}_{k,n}$ $\overset{p%
}{\rightarrow }$ $0$, therefore yield for each $k$:
\begin{equation*}
\lim_{n\rightarrow \infty }\mathcal{P}_{k,n}=\lim_{n\rightarrow \infty
}\int_{0}^{\infty }P\left( \mathcal{A}_{k,n}>\epsilon \right) d\epsilon
=\lim_{n\rightarrow \infty }\int_{0}^{1}P\left( \mathcal{A}_{k,n}>\epsilon
\right) d\epsilon =\int_{0}^{1}\lim_{n\rightarrow \infty }P\left( \mathcal{A}%
_{k,n}>\epsilon \right) d\epsilon =0.
\end{equation*}

Now apply Lemma \ref{lm:array_conv} to $\mathcal{P}_{k,n}$ to deduce that
there exists a positive integer sequence $\{\mathcal{L}_{n}\}$ that is not
unique, $\mathcal{L}_{n}$ $\rightarrow $ $\infty $ and $\mathcal{L}_{n}$ $=$
$o(n)$, such that $\lim_{n\rightarrow \infty }\mathcal{P}_{\mathcal{L}%
_{n},n} $ $=$ $\lim_{n\rightarrow \infty }\int_{0}^{1}P(\mathcal{A}_{%
\mathcal{L}_{n},n}$ $>$ $\epsilon )d\epsilon $ $=$ $0$. Therefore, by
construction $E[\mathcal{A}_{\mathcal{L}_{n},n}]$ $=$ $\int_{0}^{1}P(%
\mathcal{A}_{\mathcal{L}_{n},n}$ $>$ $\epsilon )d\epsilon $ $\rightarrow $ $%
0 $. Hence $\mathcal{A}_{\mathcal{L}_{n},n}$ $\overset{p}{\rightarrow }$ $0$
by Markov's inequality, which yields $\max_{1\leq i\leq \mathcal{L}_{n}}|%
\mathcal{X}_{n}(i)|$ $\overset{p}{\rightarrow }$ $0$ as claimed.

The sequence $\{\mathcal{L}_{n}\}$ is not unique for either of the following
reasons: ($i$) the probability limit is asymptotic hence we can always
change $\mathcal{L}_{n}$ for finitely many $n$; and ($ii$) by monotonicity
of $\max_{1\leq i\leq k}|\mathcal{X}_{n}(i)|$ any other $\{\mathcal{%
\mathring{L}}_{n}\}$ that satisfies $\mathcal{\mathring{L}}_{n}$ $%
\rightarrow $ $\infty $ and $\lim \sup_{n\rightarrow \infty }\{\mathcal{%
\mathring{L}}_{n}/\mathcal{L}_{n}\}$ $<$ $1$ satisfies $\max_{1\leq i\leq
\mathcal{\mathring{L}}_{n}}|\mathcal{X}_{n}(i)|$ $\leq $ $\max_{1\leq i\leq
\mathcal{L}_{n}}|\mathcal{X}_{n}(i)|$ $\overset{p}{\rightarrow }$ $0$ as $n$
$\rightarrow $ $\infty $.\medskip \newline
\textbf{Claim (b).}\qquad Apply the triangle inequality twice to yield both $%
\max_{1\leq i\leq k}\left\vert \mathcal{X}_{n}(i)\right\vert $ $\leq $ $%
\max_{1\leq i\leq k}|\mathcal{Y}_{n}(i)|$ $+$ $\max_{1\leq i\leq k}|\mathcal{%
X}_{n}(i)$ $-$ $\mathcal{Y}_{n}(i)|$ and $\max_{1\leq i\leq k}\left\vert
\mathcal{Y}_{n}(i)\right\vert $ $\leq $ $\max_{1\leq i\leq k}|\mathcal{X}%
_{n}(i)|$ $+$ $\max_{1\leq i\leq k}|\mathcal{X}_{n}(i)$ $-$ $\mathcal{Y}%
_{n}(i)|$, hence
\begin{equation*}
\left\vert \max_{1\leq i\leq \mathcal{L}_{n}}\left\vert \mathcal{X}%
_{n}(i)\right\vert -\max_{1\leq i\leq \mathcal{L}_{n}}\left\vert \mathcal{Y}%
_{n}(i)\right\vert \right\vert \leq \max_{1\leq i\leq \mathcal{L}%
_{n}}\left\vert \mathcal{X}_{n}(i)-\mathcal{Y}_{n}(i)\right\vert
\end{equation*}%
Now apply (a) to $\mathcal{X}_{n}(i)$ $-$ $\mathcal{Y}_{n}(i)$ to yield the
desired result. $\mathcal{QED}$.\bigskip \newline
\textbf{Proof of Theorem \ref{th:filter_X}.}\qquad By the mean value
theorem, there exists $\theta _{n}^{\ast }(i)$, $||\theta _{n}^{\ast }(i)$ $%
- $ $\theta _{0}(i)||$ $\leq $ $||\hat{\theta}_{n}(i)$ $-$ $\theta _{0}(i)||$%
, such that:
\begin{equation*}
\mathcal{X}_{n}(i,\hat{\theta}_{n}(i))=\mathcal{X}_{n}(i,\theta _{0}(i))+%
\frac{1}{\sqrt{n}}\frac{\partial }{\partial \theta ^{\prime }}\mathcal{X}%
_{n}(i,\theta _{n}^{\ast }(i))\sqrt{n}\left( \hat{\theta}_{n}(i)-\theta
_{0}(i)\right) .
\end{equation*}%
Under Assumption \ref{assum:X_filter}.d, $\sqrt{n}(\hat{\theta}_{n}(i)$ $-$ $%
\theta _{0}(i))$ $=$ $\mathcal{M}_{n}(i)$ $+$ $o_{p}(1)$ $=$ $O_{p}(1)$,
hence $||\theta _{n}^{\ast }(i)$ $-$ $\theta _{0}(i)||$ $\leq $ $||\hat{%
\theta}_{n}(i)$ $-$ $\theta _{0}(i)||$ $\overset{p}{\rightarrow }$ $0$.
Therefore $\theta _{n}^{\ast }(i)$ lies in any compact neighborhood $%
\mathcal{N}_{0}(i)$ of $\theta _{0}(i)$ with positive Lebesgue measure
asymptotically with probability approaching one. Now use Assumption \ref%
{assum:X_filter}.c and continuity of $\mathcal{D}(i,\cdot )$ to yield $%
n^{-1/2}(\partial /\partial \theta )\mathcal{X}_{n}(i,\theta _{n}^{\ast
}(i)) $ $\overset{p}{\rightarrow }$ $\mathcal{D}(i)$ $\equiv $ $\mathcal{D}%
(i,\theta _{0})$. Further, by Assumption \ref{assum:X_filter}.d $\mathcal{X}%
_{n}(i,\theta _{0}(i))$ $=$ $\mathcal{S}_{n}(i)+$ $o_{p}(1)$. We may
therefore write%
\begin{equation}
\mathcal{X}_{n}(i,\hat{\theta}_{n}(i))=\mathcal{S}_{n}(i)+\mathcal{D}%
(i)^{\prime }\mathcal{M}_{n}(i)+o_{p}(1).  \label{Xn_SnDMn}
\end{equation}

Apply Assumption \ref{assum:X_filter},d and the continuous mapping theorem
to deduce $\mathcal{S}_{n}(i)$ $+$ $\mathcal{D}(i)\mathcal{M}_{n}(i)$ $%
\overset{d}{\rightarrow }$ $\mathcal{S}(i)$ $+$ $\mathcal{D}(i)^{\prime }%
\mathcal{M}(i)$. By the definition of convergence in distribution, we may
write for each $i$:%
\begin{equation}
\mathcal{S}_{n}(i)+\mathcal{D}(i)^{\prime }\mathcal{M}_{n}(i)=\mathcal{S}(i)+%
\mathcal{D}(i)^{\prime }\mathcal{M}(i)+o_{p}(1).  \label{SnDMn_SDM}
\end{equation}%
The claim now follows from (\ref{Xn_SnDMn}) and (\ref{SnDMn_SDM}), and two
applications of Theorem \ref{th:max_p}. $\mathcal{QED}$.

\setstretch{1}
\bibliographystyle{econometrica}
\bibliography{ref_max_mean}

\end{document}